\documentclass[12pt]{amsart}
\usepackage[T1]{fontenc}
\usepackage{amssymb}
\usepackage{a4wide}
\usepackage{graphicx}
\usepackage{verbatim}
\usepackage{amsmath}
\usepackage{version}

 \theoremstyle{plain}    
 \newtheorem{thm}{Theorem}[section]
 \numberwithin{equation}{section} 
 \numberwithin{figure}{section} 
 \theoremstyle{plain}
\newtheorem{prop}[thm]{Proposition} 
 \theoremstyle{plain}    
 \newtheorem{lem}[thm]{Lemma} 
 \theoremstyle{plain}    
 \newtheorem{cor}[thm]{Corollary} 
 \theoremstyle{plain}    

 \theoremstyle{definition}
 \newtheorem{rem}[thm]{Remark} 
 \theoremstyle{definition}
 
\newcommand{\nwc}{\newcommand}

\nwc{\mf}{\mathbf} 
\nwc{\blds}{\boldsymbol} 
\nwc{\ml}{\mathcal} 

\nwc{\lam}{\lambda}
\nwc{\del}{\delta}
\nwc{\Del}{\Delta}
\nwc{\Lam}{\Lambda}
\nwc{\elll}{\ell}

\nwc{\IA}{\mathbb{A}} 
\nwc{\IB}{\mathbb{B}} 
\nwc{\IC}{\mathbb{C}} 
\nwc{\ID}{\mathbb{D}} 
\nwc{\IE}{\mathbb{E}} 
\nwc{\IF}{\mathbb{F}} 
\nwc{\IG}{\mathbb{G}} 
\nwc{\IH}{\mathbb{H}} 
\nwc{\IN}{\mathbb{N}} 
\nwc{\IP}{\mathbb{P}} 
\nwc{\IQ}{\mathbb{Q}} 
\nwc{\IR}{\mathbb{R}} 
\nwc{\IS}{\mathbb{S}} 
\nwc{\IT}{\mathbb{T}} 
\nwc{\IZ}{\mathbb{Z}} 

\def\bbleft{{\mathchoice {[\mskip-3mu {[}} {[\mskip-3mu {[}}{[\mskip-4mu {[}}{[\mskip-5mu {[}}}}
\def\bbright{{\mathchoice {]\mskip-3mu {]}} {]\mskip-3mu {]}}{]\mskip-4mu {]}}{]\mskip-5mu {]}}}}
\nwc{\setK}{\bbleft 1,K \bbright}
\nwc{\setN}{\bbleft 1,\cN \bbright}

\nwc{\va}{{\bf a}}
\nwc{\vb}{{\bf b}}
\nwc{\vc}{{\bf c}}
\nwc{\vd}{{\bf d}}
\nwc{\ve}{{\bf e}}
\nwc{\vf}{{\bf f}}
\nwc{\vg}{{\bf g}}
\nwc{\vh}{{\bf h}}
\nwc{\vi}{{\bf i}}
\nwc{\vI}{{\bf I}}
\nwc{\vj}{{\bf j}}
\nwc{\vk}{{\bf k}}
\nwc{\vl}{{\bf l}}
\nwc{\vm}{{\bf m}}
\nwc{\vM}{{\bf M}}
\nwc{\vn}{{\bf n}}
\nwc{\vo}{{\it o}}
\nwc{\vp}{{\bf p}}
\nwc{\vq}{{\bf q}}
\nwc{\vr}{{\bf r}}
\nwc{\vs}{{\bf s}}
\nwc{\vt}{{\bf t}}
\nwc{\vu}{{\bf u}}
\nwc{\vv}{{\bf v}}
\nwc{\vw}{{\bf w}}
\nwc{\vx}{{\bf x}}
\nwc{\vy}{{\bf y}}
\nwc{\vz}{{\bf z}}
\nwc{\bal}{\blds{\alpha}}
\nwc{\bep}{\blds{\epsilon}}
\nwc{\barbep}{\overline{\blds{\epsilon}}}
\nwc{\bnu}{\blds{\nu}}
\nwc{\bmu}{\blds{\mu}}
\nwc{\bet}{\blds{\eta}}

\nwc{\bk}{\blds{k}}
\nwc{\bm}{\blds{m}}
\nwc{\bM}{\blds{M}}
\nwc{\bp}{\blds{p}}
\nwc{\bq}{\blds{q}}
\nwc{\bn}{\blds{n}}
\nwc{\bv}{\blds{v}}
\nwc{\bw}{\blds{w}}
\nwc{\bx}{\blds{x}}
\nwc{\bxi}{\blds{\xi}}
\nwc{\by}{\blds{y}}
\nwc{\bz}{\blds{z}}

\nwc{\cA}{\ml{A}}
\nwc{\cB}{\ml{B}}
\nwc{\cC}{\ml{C}}
\nwc{\cD}{\ml{D}}
\nwc{\cE}{\ml{E}}
\nwc{\cF}{\ml{F}}
\nwc{\cG}{\ml{G}}
\nwc{\cH}{\ml{H}}
\nwc{\cI}{\ml{I}}
\nwc{\cJ}{\ml{J}}
\nwc{\cK}{\ml{K}}
\nwc{\cL}{\ml{L}}
\nwc{\cM}{\ml{M}}
\nwc{\cN}{\ml{N}}
\nwc{\cO}{\ml{O}}
\nwc{\cP}{\ml{P}}
\nwc{\cQ}{\ml{Q}}
\nwc{\cR}{\ml{R}}
\nwc{\cS}{\ml{S}}
\nwc{\cT}{\ml{T}}
\nwc{\cU}{\ml{U}}
\nwc{\cV}{\ml{V}}
\nwc{\cW}{\ml{W}}
\nwc{\cX}{\ml{X}}
\nwc{\cY}{\ml{Y}}
\nwc{\cZ}{\ml{Z}}

\nwc{\fA}{\mathfrak{a}}
\nwc{\fB}{\mathfrak{b}}
\nwc{\fC}{\mathfrak{c}}
\nwc{\fD}{\mathfrak{d}}
\nwc{\fE}{\mathfrak{e}}
\nwc{\fF}{\mathfrak{f}}
\nwc{\fG}{\mathfrak{g}}
\nwc{\fH}{\mathfrak{h}}
\nwc{\fI}{\mathfrak{i}}
\nwc{\fJ}{\mathfrak{j}}
\nwc{\fK}{\mathfrak{k}}
\nwc{\fL}{\mathfrak{l}}
\nwc{\fM}{\mathfrak{m}}
\nwc{\fN}{\mathfrak{n}}
\nwc{\fO}{\mathfrak{o}}
\nwc{\fP}{\mathfrak{p}}
\nwc{\fQ}{\mathfrak{q}}
\nwc{\fR}{\mathfrak{r}}
\nwc{\fS}{\mathfrak{s}}
\nwc{\fT}{\mathfrak{t}}
\nwc{\fU}{\mathfrak{u}}
\nwc{\fV}{\mathfrak{v}}
\nwc{\fW}{\mathfrak{w}}
\nwc{\fX}{\mathfrak{x}}
\nwc{\fY}{\mathfrak{y}}
\nwc{\fZ}{\mathfrak{z}}

\newcommand{\lien}{\mathfrak{n}}
\newcommand{\lienb}{\bar\lien}

\nwc{\tA}{\widetilde{A}}
\nwc{\tB}{\widetilde{B}}
\nwc{\tE}{E^{\vareps}}
\nwc{\tk}{\tilde k}
\nwc{\tN}{\tilde N}
\nwc{\tP}{\widetilde{P}}
\nwc{\tQ}{\widetilde{Q}}
\nwc{\tR}{\widetilde{R}}
\nwc{\tV}{\widetilde{V}}
\nwc{\tW}{\widetilde{W}}
\nwc{\ty}{\tilde y}
\nwc{\teta}{\tilde \eta}
\nwc{\tdelta}{\tilde \delta}
\nwc{\tlambda}{\tilde \lambda}
\nwc{\ttheta}{\tilde \theta}
\nwc{\tvartheta}{\tilde \vartheta}
\nwc{\tPhi}{\widetilde \Phi}
\nwc{\tpsi}{\tilde \psi}
\nwc{\tmu}{\tilde \mu}

\nwc{\To}{\longrightarrow} 

\nwc{\ad}{\rm ad}
\nwc{\eps}{\epsilon}
\nwc{\ep}{\epsilon}
\nwc{\vareps}{\varepsilon}

\def\ep{\epsilon}

\def\sq2{\sqrt{2}}

\def\t2{{\mathbb T}^2}
\def\s2{{\mathbb S}^2}

\nwc{\lap}{\bigtriangleup}
\nwc{\rest}{\restriction}
\nwc{\Diff}{\operatorname{Diff}}
\nwc{\diam}{\operatorname{diam}}
\nwc{\Res}{\operatorname{Res}}
\nwc{\Spec}{\operatorname{Spec}}
\nwc{\Vol}{\operatorname{Vol}}
\nwc{\Op}{\operatorname{Op}}
\nwc{\supp}{\operatorname{supp}}
\nwc{\Span}{\operatorname{span}}

\nwc{\dia}{\varepsilon}
\nwc{\cut}{f}
\nwc{\qm}{u_\hbar}

\def\hto0{\xrightarrow{\hbar\to 0}}

\def\rto0{\xrightarrow{r\to 0}}

\providecommand{\norm}[1]{\lVert#1\rVert}

\nwc{\la}{\langle}
\nwc{\ra}{\rangle}
\nwc{\lp}{\left(}
\nwc{\rp}{\right)}

\nwc{\bequ}{\begin{equation}}
\nwc{\be}{\begin{equation}}
\nwc{\ben}{\begin{equation*}}
\nwc{\bea}{\begin{eqnarray}}
\nwc{\bean}{\begin{eqnarray*}}
\nwc{\bit}{\begin{itemize}}
\nwc{\bver}{\begin{verbatim}}

\nwc{\eequ}{\end{equation}}
\nwc{\ee}{\end{equation}}
\nwc{\een}{\end{equation*}}
\nwc{\eea}{\end{eqnarray}}
\nwc{\eean}{\end{eqnarray*}}
\nwc{\eit}{\end{itemize}}
\nwc{\ever}{\end{verbatim}}

\newcommand{\bX}{\mathbf{X}}
\newcommand{\bY}{\mathbf{Y}}
\newcommand{\bS}{\mathbf{S}}

\begin{document}

\title[Products of Fourier integral operators]
{Exponential decay for products of Fourier integral operators}

\author[N. Anantharaman]{Nalini Anantharaman}
 \begin{abstract}{This text contains an alternative presentation, and in certain cases an improvement,
 of the ``hyperbolic dispersive estimate'' proved in \cite{An, AN07-1}, where it was used to make progress towards
 the quantum unique ergodicity conjecture. The main statement is a sufficient condition to have exponential decay of the norm of a product of sub-unitary Fourier integral operators. The improved version presented here is needed in the
 two papers \cite{AR10} and \cite{AS}.}\end{abstract}

\address{Laboratoire de Math\'ematique, Universit\'e d'Orsay Paris XI, 91405 Orsay Cedex}
\email{Nalini.Anantharaman@math.u-psud.fr}
    
\maketitle
  
\section{Introduction}
On a Hilbert space $\cH$, consider the product $\hat P_n\hat P_{n-1}\cdots \hat P_1$ of a large number of operators $\hat P_j$, with $\norm{\hat P_j}=1.$ Think, for instance, of the case where each operator $\hat P_j$ is an orthogonal projector, or a product of an orthogonal projector and a unitary operator. What kind of geometric considerations can be helpful to prove that the norm  $\norm{\hat P_n\hat P_{n-1}\cdots \hat P_1}$ is strictly less than $1$~? or better, that it decays exponentially fast with $n$~? In Section \ref{s:2}, we will describe a situation in which $\cH=L^2(\IR^d)$, and the operators  $\hat P_j$ are Fourier integral operators associated to a sequence of canonical transformations $\kappa_j$. We will give a ``hyperbolicity'' condition, on the sequence of transformations $\kappa_j$ and on the symbols of the operators  $\hat P_j$, under which we can prove exponential decay of the norm $\norm{\hat P_n\hat P_{n-1}\cdots \hat P_1}$.

This technique was introduced in \cite{An, AN07-1}, and is used in
\cite{An, AN07-1, AKN07, GR09, GR10, AS} to prove results related to the quantum unique ergodicity conjecture. In \cite{An, AN07-1}, the proofs are written on a riemannian manifold of negative curvature, for $\hat P_n=e^{\frac{i\tau \hbar\lap}2} \hat\chi_n$, where the operators $\hat\chi_n$ belong to a finite family of pseudodifferential operators, supported inside compact sets
of small diameters, and where $\lap$ is the laplacian and $\tau>0$ is fixed. The exponential decay is then used to prove a lower bound on the ``entropy'' of eigenfunctions, answering by the negative the long-standing question~: can a sequence of eigenfunctions concentrate on a closed geodesic, as the eigenvalue goes to infinity~? An expository paper can be found in \cite{N10}, see also the forthcoming paper \cite{ICM}. We give here an alternative presentation, based on the use of local adapted symplectic coordinates, which leads in certain cases to an improvement, needed in the two papers  \cite{AR10} and \cite{AS}.

Let us also mention the work of Nonnenmacher-Zworski \cite{NZ09, NZ09-1}, Christianson \cite{Ch1, Ch2, Ch3}, Datchev \cite{Dat}, and Burq-Guillarmou-Hassell \cite{BGH}, who showed how to use these techniques in scattering situations, to prove the existence of a gap below the real axis in the resolvent spectrum, and to get local smoothing estimates with loss, as well as Strichartz estimates. In this context, the idea of proving exponential decay for Fourier integral operators was also present, although in an implicit form, in Doi's work \cite{Doi}.

The technique is presented in the first four sections, and the applications needed in \cite{AR10, AS}
are stated in section \ref{s:exam}.

\section{A hyperbolic dispersion estimate\label{s:2}}
In this section, $\IR^d\times(\IR^d)^*$ is endowed with the canonical symplectic form $\omega_o=\sum_{j=1}^d dx_j\wedge d\xi_j $, where $dx_j$ denotes the projection on the $j$-th vector of the canonical basis in $\IR^d$, and $d\xi_j$ is the projection on the $j$-th vector of the dual basis in $(\IR^d)^*$. The space $\IR^d$ will also be endowed with its usual scalar product, denoted $\la., .\ra$, and we use it to systematically identify $\IR^d$ with $(\IR^d)^*$.

 We consider a sequence of smooth ($\cC^\infty$) canonical transformations $\kappa_n : \IR^d\times\IR^d\To \IR^d\times\IR^d$, preserving $\omega_o$ ($n\in\IN$).  We will only be interested in the restriction of $\kappa_1$ to a fixed relatively compact neighbourhood $\Omega$ of $0$, and it is actually sufficient
for us to assume that the product $\kappa_n\circ\kappa_{n-1}\circ\cdots\circ \kappa_1$ is well defined, for all $n$, on $\Omega$. The Darboux-Lie theorem ensures that every lagrangian foliation can be mapped, by a symplectic change of coordinates,
to the foliation of $\IR^d\times\IR^d$
by the ``horizontal'' leaves $\cL_{\xi_0}=\{(x, \xi) \in \IR^d\times\IR^d, \xi=\xi_0\}$. 
For our purposes (section \ref{s:exam}), there is no loss of generality if we make the simplifying assumption that each symplectic transformation $\kappa_n$ preserves this horizontal foliation. It means that $\kappa_n$ is of the form $(x, \xi)\mapsto (x', \xi'=p_n(\xi))$
 where $p_n: \IR^d\To \IR^d$ is a smooth function. In more sophisticated words, $\kappa_n$ has a generating function of the form $S_n(x, x', \theta)=
\la p_n(\theta),x'\ra-\la \theta, x\ra+\alpha_n(\theta)$ (where $x, x', \theta\in\IR^d$, and $\alpha_n:\IR^d\To \IR^d$ is a smooth function). We have the equivalence
$$\big[(x', \xi')=\kappa_n(x, \xi)\big]\Longleftrightarrow \big[\xi=-\partial_x S_n(x, x', \theta), \;\xi'=\partial_{x'} S_n(x, x', \theta),\;
\partial_\theta S_n(x, x', \theta)=0\big].$$
The product $\kappa_n\circ\ldots\circ\kappa_2\circ\kappa_1$ also preserves the horizontal foliation, and it admits the generating function \begin{multline*}\la p_n\circ\ldots\circ p_1(\theta),x'\ra-\la \theta, x\ra +\alpha_1(\theta)+\alpha_2(p_1(\theta))+\ldots+\alpha_n(p_{n-1}\circ\ldots\circ p_1(\theta))\\ =\la p_n\circ\ldots\circ p_1(\theta),x'\ra-\la \theta, x\ra +A_n(\theta),
\end{multline*}
where the equality defines $A_n(\theta)$. 

If $p$ is a map $\IR^d\To\IR^d$, we will denote $\nabla p$ the matrix $(\frac{\partial p_i}{\partial \theta_j})_{ij}$, which represents its differential in the canonical basis.

\bigskip

{\bf Assumptions (H)}~: We shall be interested in the following operators, acting on $L^2(\IR^d)$~:
$$\hat P_n f(x')=\frac1{(2\pi\hbar)^d}\int_{x\in\IR^d, \theta\in\IR^d} e^{\frac{i S_n(x, x', \theta)}\hbar}a^{(n)}(x, x', \theta, \hbar)f(x)dx d\theta,$$
where $\hbar>0$ is a parameter destined to go to $0$. We will assume the following~:
\begin{itemize}
\item[(H1)] The functions $p_n$ are smooth diffeomorphisms, and all the derivatives of $p_n$, of $p_n^{-1}$ and of $\alpha_n$ are bounded uniformly in $n$. 
\item[(H2)] For a given $\hbar>0$, the function $(x, x', \theta)\mapsto a^{(n)}(x, x', \theta, \hbar)$ is of class $\cC^\infty$;
\item[(H3)]  The function $a^{(1)}(x, x', \theta, \hbar)$ is supported in $\Omega$ with respect to the variable $x$;
\item[(H4)]  With respect to the variables $(x', \theta)$, the functions $a^{(n)}(x, x', \theta, \hbar)$ have a compact support $x'\in \Omega_1, \theta\in\Omega_2$, independent of $n$ and $\hbar$;
\item[(H5)]  When $\hbar \To 0$, each $a^{(n)}(x, x', \theta, \hbar)$ has an asymptotic expansion
$$a^{(n)}(x, x', \theta, \hbar)\sim  (\det \nabla p_n(\theta))^{1/2}\sum_{k=0}^{\infty}\hbar^k a^{(n)}_k(x, x', \theta),$$
valid up to any order and in all the $\cC^\ell$ norms on compact sets. Besides, these asymptotic expansions are {\em uniform with respect to $n$};
\item[(H6)]  If $(x', \theta')=\kappa_n(x, \theta)$, we have $|a^{(n)}_0(x, x', \theta)|\leq 1$. This condition ensures that $\norm{\hat P_n}_{L^2\To L^2}\leq 1+\cO(\hbar)$.
\end{itemize}
The operators $\hat P_n$ are (semiclassical) {\em Fourier integral operators} associated with the transformations $\kappa_n$ (section \S \ref{s:def}).

\subsection{Propagation of a single plane wave}
The following theorem is essentially proved in \cite{An}. We denote $e_{\xi_0, \hbar}$ the function $e_{\xi_0, \hbar}(x)=e^{\frac{i\la\xi_0, x\ra}\hbar}$.

\begin{thm} \label{t-main}Fix $\xi_0\in\IR^d$. Denote $\xi_n=p_n\circ\ldots\circ p_1(\xi_0)$.

In addition to the assumptions (H) above, assume that  
$$\limsup_{k\To +\infty}\frac1k\log \norm{  \nabla(p_{n+k}\circ p_{n+k-1}\circ \ldots\circ p_{n+1})(\xi_{n}) }\leq 0 ,$$
uniformly in $n$.
 
Fix $\cK>0$ arbitrary, and an integer $M\in\IN$. Then we have, for $n=\cK |\log \hbar|,$ and $\tilde\eps>0$ arbitrary,
\begin{equation*} \hat P_n\circ\ldots\circ\hat P_2\circ \hat P_1 e_{\xi_0, \hbar}(x)=e^{i\frac{ A_n(\xi_0)}\hbar}  e_{  \xi_n, \hbar}(x) (\det \nabla p_n\circ\ldots\circ p_1(\xi_0))^{1/2}
 \left[\sum_{k=0}^{M-1}\hbar^k b^{(n)}_k(x, \xi_n)\right]\\+\cO(\hbar^{M(1-\tilde\eps)}).\end{equation*}
 The functions $ b^{(n)}_k$, defined on $\IR^d\times \IR^d$,  are smooth, supported in $\Omega_1\times \Omega_2$, and
\begin{equation*}b^{(n)}_0 (x_n, \xi_n)=\prod_{j=1}^n a^{(j)}_0(x_j, x_{j+1}, \xi_j),
\end{equation*}
where we denote $\xi_n=p_n\circ\ldots\circ p_1(\xi_0)$, $x_n=x$ and the other terms are defined by the relations $(x_j, \xi_j)=\kappa_j\circ\ldots\circ\kappa_1(x_0, \xi_0)$.

The functions $b^{(n)}_k$, for $k>0$, have the same support as $b^{(n)}_0$. We have
$|b_0^{(n)}(x_n, \xi_n)|\leq 1$, and
besides, we have bounds
$$\norm{d_x^j b^{(n)}_k}_\infty\leq C(k, j,\eps)n^{j+3k}e^{\eps(j+2k)n}$$
valid for arbitrary $\eps>0$, 
where the prefactor $C(k, j, \eps)$ does not depend on $n$.

\end{thm}
If $n$ is fixed, and if we write $ \hat P_n\circ\ldots\circ\hat P_2\circ \hat P_1 e_{\xi_0, \hbar}(x)$ explicitly as an integral over
$(\IR^d)^{2n}$, this theorem is a straightforward application of the stationary phase method. If $n$ is allowed to go to infinity as $\hbar\To 0$, our result amounts, in some sense, to applying the method of stationary phase on a space whose dimension goes to $\infty$, and this is known to be very delicate. The theorem was first proved this way, in an unpublished version (available on request or on my webpage) of the paper \cite{An}. A nicer proof, written with the collaboration of St\'ephane Nonnenmacher, is available in
\cite{An}, and has also appeared under different forms in \cite{AN07-1, NZ09}. In these papers, the proofs are written on a riemannian manifold, for $\hat P_n=e^{\frac{i\tau \hbar\lap}2} \hat\chi_n$, where the operators $\hat\chi_n$ belong to a finite family of pseudodifferential operators, whose symbols are supported inside compact sets
of small diameters, and where $\lap$ is the laplacian and $\tau>0$ is fixed. In local coordinates, and on a manifold of constant negative sectional curvature, the calculations done in \cite{An, AN07-1} amount to the simpler statement presented here (see section \ref{s:exam}). 

In all the papers cited above, the dynamical systems under study satisfy a {\em uniform {hyperbolicity}} (or Anosov) property, ensuring an exponential
decay
\begin{equation}\label{e:decaynorm}\sup_{\xi\in \Omega_2} \norm{  \nabla(p_n\circ\ldots\circ p_2\circ p_1)(\xi)  }\leq C e^{-\lambda n},\end{equation}
with uniform constants $C, \lambda>0$. This is why, following \cite{NZ09}, we call our result a {\em hyperbolic dispersion estimate}.

\subsection{Estimating the norm of $\hat P_n\circ\ldots\circ\hat P_2\circ \hat P_1$}
We use the $\hbar$-Fourier transform
$$\cF_\hbar u(\xi)=\frac{1}{(2\pi\hbar)^{d/2}}\int_{\IR^d}u(x)e^{-\frac{i\la \xi, x\ra}\hbar}dx,$$
the inversion formula
$$u(x)=\frac{1}{(2\pi\hbar)^{d/2}}\int_{\IR^d}\cF_\hbar u(\xi)e^{\frac{i\la \xi, x\ra}\hbar}d\xi,$$
and the Plancherel formula $\norm{u}_{L^2(\IR^d)}=\norm{\cF_\hbar u}_{L^2(\IR^d)}.$

We denote by $\widetilde \Omega_2$ an open, relatively compact subset of $\IR^d$, that contains the closure $\overline{\Omega_2}$.
Using the Fourier inversion formula, Theorem \ref{t-main} implies, in a straightforward manner, the following
\begin{thm}\label{t:mainnorm}
In addition to the assumptions (H) above, assume that  
$$\limsup_{k\To +\infty}\frac1k\log \norm{  \nabla(p_{n+k}\circ p_{n+k-1}\circ \ldots\circ p_{n+1})(\xi_{n}) }\leq 0 ,$$
uniformly in $n$ and $\xi\in\widetilde\Omega_2$ (with $\xi_n=p_n\circ\ldots\circ p_1(\xi)$).
 
 Fix $\cK>0$ arbitrary. Then, for $n=\cK |\log \hbar|$,  
$$\norm{ \hat P_n\circ\ldots\circ\hat P_2\circ \hat P_1}_{L^2\To L^2}\leq {\frac{|\widetilde\Omega_2|^{1/2}}{(2\pi\hbar)^{d/2}}}
\sup_{\xi\in\widetilde\Omega_2} |\det \nabla p_n\circ\ldots\circ p_1(\xi)|^{1/2} (1+\cO(\hbar n^3 e^{\eps n})),$$
 where $|\widetilde\Omega_2|$ denotes the volume of $\widetilde\Omega_2$ (and $\eps>0$ is arbitrary).
\end{thm}
Of course, since $\norm{\hat P_j}_{L^2\To L^2}\leq 1+\cO(\hbar)$, we always have the trivial bound $\norm{ \hat P_n\circ\ldots\circ\hat P_2\circ \hat P_1}_{L^2\To L^2}\leq 1+\cO(\hbar|\log\hbar|)$.
Since we are working in the limit $\hbar \To 0$, our estimate can only have an interest if we have an upper bound of the form
\begin{equation}\label{e:decaydet}\sup_{\xi\in\widetilde\Omega_2} |\det \nabla p_n\circ\ldots\circ p_1(\xi)|^{1/2}\leq C e^{-\lambda n}, \qquad \lambda>0,\end{equation}
and if $\cK$ is large enough. Note that \eqref{e:decaydet} is weaker than the condition
\eqref{e:decaynorm}.

We now state a refinement of Theorem \ref{t:mainnorm}. We consider the same family $\hat P_i$, satisfying assumptions (H). The multiplicative constants in our estimate have no importance, and in what follows we will often omit them. 

\begin{thm}\label{t:mainnorm2} 
 Assume as above that
$$\limsup_{k\To +\infty}\frac1k\log \norm{  \nabla(p_{n+k}\circ p_{n+k-1}\circ \ldots\circ p_{n+1})(\xi_{n}) }\leq 0 ,$$
uniformly in $n$ and $\xi\in\widetilde \Omega_2$.

Let $r\leq d$, and assume that the coisotropic foliation by the leaves $\{ \xi_{r+1}=c_{r+1}, \ldots, \xi_d=c_d\}$ is invariant by each canonical transformation
$\kappa_n$. In other words, the map $p_n$ is of the form
$$p_n((\xi_1,\ldots, \xi_r), (\xi_{r+1}, \ldots, \xi_d))  =\left(m_n (\xi_1,\ldots, \xi_d), 
 \tilde p_n(\xi_{r+1}, \ldots, \xi_d)\right),$$
 where $m_n:\IR^d\To \IR^{r}$ and $\tilde p_n:\IR^{d-r}\To \IR^{d-r}.$

Fix $\cK>0$ arbitrary. Then there exists $\hbar_\cK>0$ such that, for $n=\cK |\log \hbar|,$ and for $\hbar<\hbar_\cK$,
$$\norm{ \hat P_n\circ\ldots\circ\hat P_2\circ \hat P_1}_{L^2\To L^2}\leq {\frac{1}{(2\pi\hbar)^{(r-\eps)/2}}}
 \frac{ \sup_{\xi\in\widetilde \Omega_2}|(\det \nabla  p_n\circ\ldots\circ p_1(\xi))|^{1/2}}{
\inf_{\xi\in\widetilde \Omega_2}|(\det \nabla\tilde    p_n\circ\ldots\circ \tilde p_1(\xi))|^{1/2}} (1+\cO(n^3\hbar e^{\eps n}))$$
for $\eps>0$ arbitrary.

In addition, if we make the stronger assumption that $ \norm{  \nabla(p_{n+k}\circ p_{n+k-1}\circ \ldots\circ p_{n+1})(\xi_{n}) } $ is bounded above, 
uniformly in $n, k$ and for $\xi\in\tilde\Omega_2$, we have
$$\norm{ \hat P_n\circ\ldots\circ\hat P_2\circ \hat P_1}_{L^2\To L^2}\leq {\frac{1}{(2\pi\hbar)^{r/2}}}
 \frac{ \sup_{\xi\in\widetilde \Omega_2}|(\det \nabla  p_n\circ\ldots\circ p_1(\xi))|^{1/2}}{
\inf_{\xi\in\widetilde \Omega_2}|(\det \nabla\tilde    p_n\circ\ldots\circ \tilde p_1(\xi))|^{1/2}} (1+\cO(n^3\hbar)).$$

 \end{thm}
 Theorem \ref{t:mainnorm2} is an improvement of Theorem \ref{t:mainnorm} in the case where
 we have
$$\frac{ 1}{(2\pi\hbar)^{d/2}}
\sup_{\xi\in\Omega_2} |(\det \nabla p_n\circ\ldots\circ p_1(\xi_0))^{1/2}|\gg 1$$ but 
$${\frac{1}{(2\pi\hbar)^{r/2}}}
 \frac{ \sup_{\xi\in\Omega_2}|(\det \nabla  p_n\circ\ldots\circ p_1(\xi))|^{1/2}}{
\inf_{\xi\in\Omega_2}|(\det \nabla\tilde    p_n\circ\ldots\circ \tilde p_1(\xi))|^{1/2}}\ll1.$$ As a trivial
example, when each $\kappa_n$ is the identity, Theorem \ref{t:mainnorm}
gives a non-optimal bound, whereas we can take $r=0$ in Theorem \ref{t:mainnorm2}, and recover
the (almost) optimal bound $\norm{ \hat P_n\circ\ldots\circ\hat P_2\circ \hat P_1}_{L^2\To L^2}\leq 1+\cO(\hbar|\log\hbar|^3).$ A less trivial example will be given in section \ref{s:exam}.

\section{Proof of Theorem \ref{t-main}}
The ideas below are contained in \cite{An, AN07-1}; however, our notations here are quite
different, and we recall (without giving all details) the main steps. In all this section, $M$ is a fixed
integer, and all the calculations are done modulo remainders of order $\hbar^{M}$
(with explicit control of the constants). 

In all that follows, it is useful to keep in mind the following~: if $(x', \xi')=\kappa_n\circ\ldots\circ\kappa_2\circ\kappa_1(x, \xi)$, we have $\xi'=p_n\circ\ldots\circ p_2\circ p_1(\xi)$, and $x=  \nabla (p_n\circ\ldots\circ p_2\circ p_1)^\intercal x'+\nabla A_n(\xi)$.

\subsection{One step of the iteration\label{s:onestep}}
Let us first fix $\xi\in\IR^d$, and look at the action of the operator $\hat P_n$ on a function of the form
$$b_{\xi}(x)=e^{\frac{i\la \xi, x\ra}\hbar}b(x)$$
where $$b(x)=\sum_{k=0}^{M-1}\hbar^k b_k(x),$$
and where the functions $b_k$ are of class $\cC^\infty$.

We introduce the following notation~: 
$$(T_n^\xi a)(x')=a^{(n)}_0(x, x', \xi)a(x)$$
where $x$ is the point such that $(x',  p_n(\xi))=\kappa_n(x, \xi)$ (in other words, $x= \nabla p_n(\xi)^\intercal x'+\nabla \alpha_n(\xi)$). In the case $a^{(n)}_0\equiv 1$, we note that the operator $U_n^\xi: a\mapsto (\det \nabla p_n(\xi))^{1/2}T_n^\xi a$ is unitary on $L^2(\IR^d)$. If we assume (as above) that $|a^{(n)}_0(x, x', \xi)|\leq 1$, it defines a bounded operator on $L^2(\IR^d)$,
of norm $\leq 1$.

A standard application of the stationary phase method yields~:
\begin{prop}\label{p:onestep}
$$\hat P_nb_{\xi}(x')=
e^{i\frac{ \alpha_n(\xi)+ \la p_n(\xi), x'\ra}\hbar}   (\det \nabla p_n(\xi))^{1/2}
 \left[\sum_{k=0}^{M-1}\hbar^k b'_k(x')\right]+\hbar^{M}R_{M}(x'),$$
 where~:
 \begin{itemize}
 \item $b'_0(x')=(T_n^\xi b_0)(x')$;
 \item $b'_k(x')=\sum_{0\leq l \leq k-1} D^{2(k-l)}_n b_l(x')+(T_n^\xi b_k)(x')$, where the operator
 $ D^{2(k-l)}_n$ is a differential operator of order $2(k-l)$ (whose expression also depends on $\xi$, although it does not appear in our notations). Its coefficients can be expressed in terms of the derivatives of order $\leq 2(k-l)$ of $a_l^{(n)}$, and of order $\leq 2(k-l)+3$ of $p_n, p_n^{-1}$ and $\alpha_n$, at the point $(x, x', \xi)$, where $(x',  p_n(\xi))=\kappa_n(x, \xi)$).
 \item There exists an integer $N_d$ (depending only on the dimension $d$), and a positive
 real number $C$ such that
 $$\norm{R_{M}}_{L^2(\IR^d)}\leq C \sum_{k=0}^{M-1} \norm{b_k}_{\cC^{2(M-k)+N_d}}.$$
 \end{itemize}

\end{prop} The constant $C$ can be expressed in terms of a fixed finite number of derivatives of the functions $a_l^{(n)}$ ($l\leq M-1$), $p_n, p_n^{-1}$ and $\alpha_n$ at the point $(x, x', \xi)$.
 Under our assumptions (H1) and (H5), $C$ is uniformly bounded for all $n$. Also note that, under (H4), the functions $b'_k$ are always supported inside the relatively compact set $\Omega_1$.

\subsection{After many iterations}
We can now describe the action  of the product $ \hat P_n\circ\ldots\circ\hat P_2\circ \hat P_1$ on
$e_{\xi_0, \hbar}$. We will give an approximate expression of $ \hat P_n\circ\ldots\circ\hat P_2\circ \hat P_1 e_{\xi_0, \hbar}(x)$, in the form
$$e^{i\frac{ A_n(\xi_0)}\hbar}  e_{  \xi_n, \hbar}(x) (\det \nabla p_n\circ\ldots\circ p_1(\xi_0))^{1/2}
 \left[\sum_{k=0}^{M-1}\hbar^k b^{(n)}_k(x)\right],$$
as announced in the theorem. This expression will approximate $ \hat P_n\circ\ldots\circ\hat P_2\circ \hat P_1 e_{\xi_0, \hbar}$
up to an error of order $\hbar^{M(1-\tilde \eps)}$ for any $\tilde\eps>0$. The function $b^{(n)}_k(x)$ depends, of course, on $\xi_0$, and in the final statement of the theorem we indicated this dependence by writing $b^{(n)}_k(x,\xi_n)$ (with $\xi_n= p_n\circ\ldots\circ p_1(\xi_0)$).

The method consists in iterating the method described in Section \ref{s:onestep}, controlling
carefully how the remainders grow with $n$ in the $L^2$ norm. We recall that $\norm{\hat P_n}_{L^2(\IR^d)}
\leq 1+\cO(\hbar)$, uniformly in $n$. 

Suppose that, after $n$ iterations, we have proved that
$$\hat P_n\circ\ldots\circ\hat P_2\circ \hat P_1 e_{\xi_0, \hbar}(x)=e^{i\frac{ A_n(\xi_0)}\hbar}  e_{  \xi_n, \hbar}(x) (\det \nabla p_n\circ\ldots\circ p_1(\xi_0))^{1/2}
 \left[\sum_{k=0}^{M-1}\hbar^k b^{(n)}_k(x)\right]+\hbar^{M}\cR^{(n)}_M(x).$$

  The calculations done in Section \ref{s:onestep} allows to describe the action
 of $\hat P_{n+1}$ on $e_{  \xi_n, \hbar}  
 \left[\sum_{k=0}^{M-1}\hbar^k b^{(n)}_k\right]$~:
 \begin{multline}\label{e:iter}\hat P_{n+1}e_{  \xi_n, \hbar}  
 \left[\sum_{k=0}^{M-1}\hbar^k b^{(n)}_k\right](x)\\
 =
 e^{i\frac{ \alpha_{n+1}(\xi_n)+ \la p_{n+1}(\xi_n), x\ra}\hbar}   (\det \nabla p_{n+1}(\xi_n))^{1/2}
 \left[\sum_{k=0}^{M-1}\hbar^k b^{(n+1)}_k(x)\right]+\hbar^{M}R^{(n+1)}_{M}(x).
 \end{multline}
 Note that
$$ \prod_{\ell=1}^n(\det \nabla p_{\ell}(\xi_{\ell-1}))^{1/2}= (\det \nabla p_n\circ\ldots\circ p_1(\xi_0))^{1/2},$$
and
$$A_n(\xi_0)=\alpha_1(\xi_0)+\alpha_2(\xi_1)+\cdots+\alpha_n(\xi_{n-1}),$$
so that 
\begin{multline*}\hat P_{n+1}\hat P_n\circ\ldots\circ\hat P_2\circ \hat P_1 e_{\xi_0, \hbar}(x)=e^{i\frac{ A_{n+1}(\xi_0)}\hbar}  e_{  \xi_{n+1}, \hbar}(x) (\det \nabla p_{n+1}\circ\ldots\circ p_1(\xi_0))^{1/2}
 \left[\sum_{k=0}^{M-1}\hbar^k b^{(n+1)}_k(x)\right]\\+\hbar^{M}\cR^{(n+1)}_M(x),\end{multline*}
 with the relation $\cR_M^{(n+1)}=e^{i\frac{ A_n(\xi_0)}\hbar} (\det \nabla p_n\circ\ldots\circ p_1(\xi_0))^{1/2}R^{(n+1)}_M+\hat P_{n+1}\cR^{(n)}_M$.

 We need to control how each term in these expansions will grow with $n$, and in particular, to control the remainder terms.
  We form an array $B^{(n)}$ that contains all the functions $b^{(n)}_k$, and a certain number of higher order differentials~:
 $$B^{(n)}_{j, k}=d^jb^{(n)}_k,$$
 with $0\leq k\leq M-1$ and $0\leq j\leq 2(M-k)+N_d.$ The index $k$ indicates the power of $\hbar$, and the index $j$ indicates the number of differentials. Note that $d^jb^{(n)}_k$ is a (symmetric) covariant tensor field of order $j$ on $\IR^d$.  If $\sigma$ is a covariant tensor field of order $j$ on $\IR^d$, we define $\norm{\sigma}_\infty=\sup_{x\in \IR^d}|\sigma_x|,$ where $|\sigma_x|$ is the norm of the $j$-linear form
 $\sigma_x$.
 By assumption (H4), the forms $d^jb^{(n)}_k$ all vanish outside the compact set $\Omega_1$.
 
  There is a linear relation between $B^{(n)}$ and $B^{(n+1)}$, that we can make a little more explicit. 
We extend the definition of the operators $T_n^\xi$ (previously defined on functions) to covariant tensor fields, by letting (if $\sigma$ is of order $j$)
$$(T_n^\xi\sigma)_{x'}(v_1,\ldots, v_j)= 
a^{(n)}_0(x, x', \xi)\sigma_{x}( \nabla p_n(\xi)^\intercal v_1,\ldots, \nabla p_n(\xi)^\intercal v_j),$$
where $x= \nabla p_n(\xi)^\intercal x'$. Taking successive derivatives of the relation
$$b^{(n+1)}_k=\sum_{0\leq l \leq k-1} D^{2(k-l)}_{n+1} b^{(n)}_l+T_{n+1}^\xi b_k^{(n)},$$
which appears in Proposition \ref{p:onestep},
we obtain a linear relation of the form~:
$$B^{(n+1)}=  K_{n+1} B^{(n)}+ L_{n+1}  B^{(n)}
+ T_{n+1}^\xi B^{(n)}, $$
where  $T_{n+1}^\xi $ acts ``diagonally'', meaning that $[T_{n+1}^\xi B^{(n)}]_{j, k}=T_{n+1}^\xi \left(B^{(n)}_{j, k}\right).$ The only information we need about the other terms is that
$[K_{n+1} B^{(n)}]_{j, k}$ depends only on the components  $B^{(n)}_{j', l}$,
for $l\leq k-1$ and $j'\leq 2(k-l)+j$; and $[L_{n+1}B^{(n)}]_{j, k}$ depends only on the components  $B^{(n)}_{j', k}$, with $j'\leq j-1$. Besides, we have
$$\max_{j, k}\norm{[K_{n+1} B^{(n)}]_{j, k}}_\infty\leq C \max_{j, k}\norm{  B^{(n)}_{j, k}}_\infty,$$
where $C$ does not depend on $n$ by our assumptions (H4) (and the same holds with $K_{n+1}$ replaced by $L_{n+1}$).

By induction, we see that $B^{(n)}$ can be expressed as
$$B^{(n)}=\sum_{A_{\ell}\in\{ T_\ell^\xi, K_\ell, L_\ell\} }A_{n}\circ A_{n-1}\circ \cdots\circ A_{1} B^{(0)}.$$
In a product of the form $A_{n}\circ A_{n-1}\circ \cdots\circ A_{1} $ (where $A_{\ell}\in\{ T_\ell^\xi, K_\ell, L_\ell\} $ for all $\ell=1,\ldots, n$), we see that there can be at most $M$ indices $\ell$ for which $A_\ell=K_\ell$, and $2M+ N_d$ indices $k$ such that  $A_\ell= L_\ell$ (otherwise the product $A_{n}\circ A_{n-1}\circ \cdots\circ A_{1} $ vanishes). Even more precisely, when we write
\begin{equation}\label{e:iterate}B^{(n)}_{j, k}=\left[\sum_{A_{\ell}\in\{ T_\ell^\xi, K_\ell, L_\ell\} }A_{n}\circ A_{n-1}\circ \cdots\circ A_{1} B^{(0)}\right]_{j, k},
\end{equation}
in the right-hand side there can be at most $k$ indices $\ell$ with $A_\ell=K_\ell$, and $2k+j$ indices
$\ell$ with $A_\ell=L_\ell$. Hence, the sum has at most $2^{3k+j}C_n^{3k+j}\sim C(k, j) n^{3k+j}$ terms.

 We now use our assumption that
$$\limsup_{k\To +\infty}\frac1k\log \norm{  \nabla(p_{n+k}\circ p_{n+k-1}\circ \ldots\circ p_{n+1})(\xi_{n}) }\leq 0 ,$$
uniformly in $n$. Combined with \eqref{e:iterate}, this implies that, for any $\eps>0$, we have
$$\norm{B^{(n)}_{j, k}}_\infty \leq C(k, j, \eps) n^{3k+j} e^{\eps n (j+2k)}.$$
These estimates (combined with Proposition \ref{p:onestep}) imply that 
$$\norm{R^{(n+1)}_M}_{L^2(\IR^d)}\leq C\sum_{k=0}^{M-1}\sum_{j=0}^{2(M-k)+N_d}\norm{B^{(n)}_{j, k}}_\infty\leq C(M, j, \eps)n^{3M+N_d}e^{\eps n(2M+N_d)}.$$
Remember the induction relation
$$\cR_M^{(n+1)}=e^{i\frac{ A_n(\xi_0)}\hbar} (\det \nabla p_n\circ\ldots\circ p_1(\xi_0))^{1/2}R^{(n+1)}_M+\hat P_{n+1}\cR^{(n)}_M.$$ 
We have
 $\norm{\hat P_{n+1}\cR^{(n)}_M}_{L^2(\IR^d)}\leq (1+\cO(\hbar))\norm{\cR^{(n)}_M}_{L^2(\IR^d)}.$  
If we restrict our attention to $n\leq \cK|\log\hbar|$ (where $\cK$ is fixed), this induction relation implies that $\norm{\cR^{(n)}_M}_{L^2(\IR^d)}=\cO(\hbar^{-\tilde\eps M})$ for any $\tilde\eps>0$.


\section{Proof of Theorems \ref{t:mainnorm} and \ref{t:mainnorm2}}
\subsection{Theorem \ref{t:mainnorm}}
The proof of Theorem \ref{t:mainnorm} is now very easy. Let $u\in L^2(\IR^d)$. We know that
$$u(x)=\frac{1}{(2\pi\hbar)^{d/2}}\int_{\IR^d}\cF_\hbar u(\xi)e^{\frac{i\la \xi, x\ra}\hbar}d\xi.$$
Let $\widetilde\Omega_2$ be an open set containing the closure of $\Omega_2$. We decompose
$u=u_1+u_2$, where
$$u_1(x)=\frac{1}{(2\pi\hbar)^{d/2}}\int_{\widetilde\Omega_2}\cF_\hbar u(\xi)e^{\frac{i\la \xi, x\ra}\hbar}d\xi$$
 and
 $$u_2(x)=\frac{1}{(2\pi\hbar)^{d/2}}\int_{\IR^d\setminus \widetilde\Omega_2}\cF_\hbar u(\xi)e^{\frac{i\la \xi, x\ra}\hbar}d\xi.$$
 Since $\hat P_1^*\hat P_1$ is a pseudodifferential operator, whose complete symbol is supported
 in $\Omega\times \Omega_2$, we have $\norm{\hat P_1 u_2}_{L^2(\IR^d)}=\cO(\hbar^\infty)\norm{u_2}_{L^2(\IR^d)}.$
 
 Concerning $u_1$, we apply Theorem \ref{t-main} for each $\xi\in \widetilde\Omega_2$. We take $n=\cK|\log\hbar|$ and choose $M$ accordingly, large enough so that
 $$\cO(\hbar^{M(1-\tilde \eps)})\ll  \sup_{\xi\in\widetilde\Omega_2} |\det \nabla p_n\circ\ldots\circ p_1(\xi)|^{1/2}.$$
  From Theorem \ref{t-main}, we know that
 $$\norm{ \hat P_n\circ\ldots\circ\hat P_2\circ \hat P_1 e_{\xi, \hbar}}_{L^2(\IR^d)}\leq 
 |\det \nabla p_n\circ\ldots\circ p_1(\xi)|^{1/2} (1+\cO(\hbar n^3 e^{2\eps n}))$$
(for $\eps>0$ arbitrary). By a direct application of the triangular inequality, it
 follows that
 \begin{multline*}\norm{ \hat P_n\circ\ldots\circ\hat P_2\circ \hat P_1 u_1}_{L^2(\IR^d)}\leq 
 \frac{1}{(2\pi\hbar)^{d/2}}
  \sup_{\xi\in\widetilde\Omega_2} |\det \nabla p_n\circ\ldots\circ p_1(\xi)|^{1/2} (1+\cO(\hbar n^3 e^{2\eps n}))\norm{\cF_\hbar u}_{L^1(\widetilde\Omega_2)}\\
  \leq \frac{1}{(2\pi\hbar)^{d/2}}
  \sup_{\xi\in\widetilde\Omega_2} |\det \nabla p_n\circ\ldots\circ p_1(\xi)|^{1/2} (1+\cO(\hbar n^3 e^{2\eps n}))|\widetilde\Omega_2|^{1/2}\norm{\cF_\hbar u}_{L^2(\IR^d)}
  \end{multline*}
  and our result follows.
  
  \subsection{Theorem \ref{t:mainnorm2}\label{s:proof}}
\subsubsection{The Cotlar-Stein lemma}   
  \begin{lem}Let $E, F$ be two Hilbert spaces. Let $(A_\alpha)\in \cL(E, F)$ be a countable family
of bounded linear operators from $E$ to $F$. Assume that for some $R>0$ we have
$$\sup_{\alpha}\sum_{\beta}\norm{A^*_\alpha A_\beta}^{\frac12}\leq R$$ and
$$\sup_{\alpha}\sum_{\beta}\norm{A_\alpha A^*_\beta}^{\frac12}\leq R$$
Then $A=\sum_{\alpha}A_\alpha$ converges strongly and $A$ is a bounded operator with $\norm{A}\leq R$.
\end{lem}
The Cotlar-Stein lemma is often used to bound in a precise manner the norm of pseudodifferential operators.

\subsubsection{\label{s:decompo}}  Remember that we assume everywhere that $n=\cK|\log\hbar|$, with $\cK$ fixed. In order to bound the norm of
$\hat P_n\circ\ldots\circ\hat P_1$ (modulo $\hbar^N$ for arbitrary $N$), the results of the previous sections show that it is enough to bound the norm of the operator $A$ defined by
\begin{multline}\label{e:initial}Af(x')\\=\frac{1}{(2\pi\hbar)^d}\int_{\xi\in\widetilde\Omega_2, x\in\IR^d} 
 (\det \nabla p_n\circ\ldots\circ p_1(\xi))^{1/2}
 \left[\sum_{k=0}^{M-1}\hbar^k b^{(n)}_k(x', \xi_n)\right]e^{\frac{i}\hbar\left(\la \xi_n, x'\ra+A_n(\xi)-\la \xi, x\ra\right)} f(x)dxd\xi\\
 =\frac{1}{(2\pi\hbar)^{d/2}}\int_{\xi\in\widetilde\Omega_2, x\in\IR^d} 
 (\det \nabla p_n\circ\ldots\circ p_1(\xi))^{1/2}
 \left[\sum_{k=0}^{M-1}\hbar^k b^{(n)}_k(x', \xi_n)\right]e^{\frac{i}\hbar\left(\la \xi_n, x'\ra+A_n(\xi)\right)} \cF_\hbar f(\xi)d\xi
 \end{multline}
 for a suitable choice of $M$, large.
 We denote everywhere $\xi_n=p_n\circ\ldots\circ  p_1(\xi).$
 
 We decompose $\IR^d=\IR^r\times \IR^{d-r}$, and write any $\xi\in \IR^d$ as $\xi=(\xi_{(r)}, \tilde\xi)$ where $\xi_{(r)}\in \IR^r$ and $ \tilde\xi\in  \IR^{d-r}$. Under our current assumptions, $\xi_n$ decomposes as
  $\xi_n=(\xi_{n(r)}, \tilde\xi_n)$, where $\tilde\xi_n=\tilde p_n\circ\ldots\circ \tilde p_1(\tilde\xi).$
  
  We now introduce a (real-valued) smooth compactly supported $\chi$ on $ \IR^{d-r}$, such that $0\leq\chi\leq 1$, and having the property that
  $$\sum_{\ell\in \IZ^{d-r}}\chi(\tilde\xi-\ell)=1$$
  for all $\tilde\xi\in  \IR^{d-r}$. For $\hbar>0$, $\ell\in \IZ^{d-r}$ and $\tilde\xi\in  \IR^{d-r}$, we denote $\chi_{\hbar, \ell}(\tilde\xi)=\chi\left(\frac{\tilde\xi}{2\pi\hbar}-\ell\right)$. Using the same notation as in \eqref{e:initial}, we define
  \begin{multline}\label{e:decompo}A_\ell f(x')\\=\frac{1}{(2\pi\hbar)^d}\int 
 (\det \nabla p_n\circ\ldots\circ p_1(\xi))^{1/2}
 \left[\sum_{k=0}^{M-1}\hbar^k b^{(n)}_k(x', \xi_n)\right]\chi_{\hbar, \ell}(\tilde\xi_n)e^{\frac{i}\hbar\left(\la \xi_n, x'\ra+A_n(\xi)-\la \xi, x\ra\right)} f(x)dxd\xi\\
 =\frac{1}{(2\pi\hbar)^{d/2}}\int_{\xi\in\widetilde\Omega_2}  
 (\det \nabla p_n\circ\ldots\circ p_1(\xi))^{1/2}
 \left[\sum_{k=0}^{M-1}\hbar^k b^{(n)}_k(x', \xi_n)\right]\chi_{\hbar, \ell}(\tilde\xi_n)e^{\frac{i}\hbar\left(\la \xi_n, x'\ra+A_n(\xi)\right)} \cF_\hbar f(\xi)d\xi.
 \end{multline}
 It is clear that $A=\sum_{\ell\in \IZ^{d-r}}A_\ell.$ A crucial remark is that the function $\xi\mapsto \chi_{\hbar, \ell}(\tilde\xi_n)$, defined on $\Omega_2$, is supported in a set of volume $\leq (2\pi\hbar)^{d-r}
\frac{1}{ \inf_{\xi\in\widetilde \Omega_2}|(\det \nabla\tilde    p_n\circ\ldots\circ \tilde p_1(\xi))|}.$

 We are going to apply the Cotlar-Stein lemma to this decomposition. Let us write explicitly the expression for the adjoint~:
  \begin{multline}\label{e:adjoint}A^*_\ell f(x)\\=\frac{1}{(2\pi\hbar)^d}\int 
 (\det \nabla p_n\circ\ldots\circ p_1(\xi))^{1/2}
 \left[\sum_{k=0}^{M-1}\hbar^k\overline{ b^{(n)}_k}(x', \xi_n)\right]\chi_{\hbar, \ell}(\tilde\xi_n)e^{-\frac{i}\hbar\left(\la \xi_n, x'\ra+A_n(\xi)-\la \xi, x\ra\right)} f(x')dx'd\xi .
 \end{multline}
We shall evaluate the norm of $A^*_mA_\ell$ and $A_\ell A^*_m$, for all $m, \ell\in \IZ^{d-r}$.

\subsubsection{Norm of $A^*_mA_\ell$}
We evaluate the norm of $A^*_mA_\ell$ acting on $L^2(\IR^d)$ by studying the scalar product
$\la A_\ell f, A_m f\ra$ for $f\in L^2(\IR^d)$. Using expression \eqref{e:decompo} and bilinearity of the scalar product, we will bound the scalar product $\la A_\ell f, A_m f\ra$ by studying separately each bracket
\begin{equation}\label{e:scalar}\chi_{\hbar, \ell}(\tilde\xi_n)\chi_{\hbar, m}(\tilde\xi'_n)\left\la  \left[\sum_{k=0}^{M-1}\hbar^k b^{(n)}_k(x', \xi_n)\right]e^{\frac{i}\hbar\la \xi_n, x'\ra },
 \left[\sum_{k=0}^{M-1}\hbar^k b^{(n)}_k(x', \xi'_n)\right]e^{\frac{i}\hbar\la \xi'_n, x'\ra }
 \right\ra_{L^2_{x'}}.\end{equation}
 Using the notation of \S  \ref{s:decompo}, we decompose the complex phase $\la \xi_n, x'\ra -\la \xi'_n, x'\ra$ into $ \la \xi_{n(r)}, x'_{(r)}\ra -\la \xi'_{n (r)}, x'_{(r)}\ra + \la \tilde\xi_n, \tilde x'\ra -\la \tilde\xi'_n, \tilde x'\ra.$ In the integral defining the scalar product \eqref{e:scalar}, we perform an integration by parts with respect to $\tilde x'\in\IR^{d-r}$~: we integrate $N$ times the function $e^{\frac{i}\hbar \la \tilde\xi_n, \tilde x'\ra -\la \tilde\xi'_n, \tilde x'\ra}$ and differentiate the functions $b^{(n)}_k(x', \xi_n)$. Using the estimates of Theorem \ref{t-main}, we obtain
  \begin{prop}
 \begin{multline}\label{e:scalarbound}\chi_{\hbar, \ell}(\tilde\xi_n)\chi_{\hbar, m}(\tilde\xi'_n)\left\lvert\left\la  \left[\sum_{k=0}^{M-1}\hbar^k b^{(n)}_k(x', \xi_n)\right]e^{\frac{i}\hbar\la \xi_n, x'\ra },
 \left[\sum_{k=0}^{M-1}\hbar^k b^{(n)}_k(x', \xi'_n)\right]e^{\frac{i}\hbar\la \xi'_n, x'\ra }
 \right\ra \right\rvert \\
 \leq C(\eps)e^{\eps Nn}\frac{1}{(\norm{m-\ell}+1)^N} 
  \end{multline}
  for $\eps>0$ arbitrary.
 \end{prop}
 The integer $N$ will be chosen soon.
We now use the bilinearity of the scalar product, and the fact that
$$\norm{\chi_{\hbar, \ell}(\tilde\xi_n)\cF_\hbar f(\xi)}_{L^1(\widetilde\Omega_2)}\leq 
(2\pi\hbar)^{(d-r)/2}
 \frac{ 1}{
\inf_{\xi\in\widetilde\Omega_2}|(\det \nabla\tilde    p_n\circ\ldots\circ \tilde p_1(\xi))|^{1/2}} 
\norm{\cF_\hbar f(\xi)}_{L^2(\IR^d)}.$$ Combined with expression \eqref{e:decompo}, this yields that
\begin{equation}\label{e:norm1}\norm{A^*_mA_\ell}\leq C(\eps)e^{\eps Nn}
\frac{1}{(\norm{m-\ell}+1)^N}  {\frac{1}{(2\pi\hbar)^{r}}}
 \frac{ \sup_{\xi\in\widetilde\Omega_2}|(\det \nabla  p_n\circ\ldots\circ p_1(\xi))|}{
\inf_{\xi\in\widetilde\Omega_2}|(\det \nabla\tilde    p_n\circ\ldots\circ \tilde p_1(\xi))|} .
\end{equation}
Looking at the statement of the Cotlar-Stein lemma, we see that we must choose $N$ large enough such that $\sum_{\ell\in\IZ^{d-r}}\frac{1}{(\norm{\ell}+1)^{N/2}}<+\infty.$

\begin{rem}\label{r:unif}If we make the assumption that $ \norm{  \nabla(p_{n+k}\circ p_{n+k-1}\circ \ldots\circ p_{n+1})(\xi_{n}) } $ is bounded above, 
uniformly in $n, k$ and $\xi\in\widetilde\Omega_2$, we see that we can take $\eps=0$ in all the statements made above.
\end{rem}

\subsubsection{Norm of $A_\ell A^*_m$}
This step is actually shorter than the previous one. We now have to evaluate the scalar product
$\la A^*_\ell f, A^*_m f\ra$ for $f\in L^2(\IR^d)$, and we use the expression \eqref{e:adjoint} of the adjoint.  We do not need integration by parts, as we see directly that $\la A^*_\ell f, A^*_m f\ra$ vanishes as soon as $\norm{m-\ell}$ is too large (in fact, the supports of $\chi_{\hbar, \ell}$ and $\chi_{\hbar, m}$ are disjoint if $\norm{m-\ell}>C$, where $C$ is fixed and depends only on the support of $\chi$). In what follows we consider the case $\norm{m-\ell}\leq C$. We see that $A^*_\ell f$ is the $\cF_\hbar$-transform of
$$F_\ell: \xi\mapsto 
\frac{1}{(2\pi\hbar)^{d/2}}\int 
 (\det \nabla p_n\circ\ldots\circ p_1(\xi))^{1/2}
 \left[\sum_{k=0}^{M-1}\hbar^k\overline{ b^{(n)}_k}(x', \xi_n)\right]\chi_{\hbar, \ell}(\tilde\xi_n)e^{-\frac{i}\hbar\left(\la \xi_n, x'\ra+A_n(\xi)\right)} f(x')dx'.$$
 We recall that each $b^{(n)}_k(x', \xi_n)$ is supported in $\{x'\in\Omega_1\}$, and we bound
 \begin{multline*}\norm{F_\ell}_{L^2(\IR^d)}\leq \\
 \frac{1}{(2\pi\hbar)^{d/2}}\sup_{\xi\in\widetilde\Omega_2}|(\det \nabla  p_n\circ\ldots\circ p_1(\xi))|^{1/2}
 (2\pi\hbar)^{(d-r)/2} \frac{1}{
\inf_{\xi\in\widetilde\Omega_2}|(\det \nabla\tilde    p_n\circ\ldots\circ \tilde p_1(\xi))|^{1/2}}\norm{f}_{L^1(\Omega_1)},
\end{multline*}
and $\norm{f}_{L^1(\Omega_1)}\leq |\Omega_1|^{1/2}\norm{f}_{L^2(\IR^d)}.$
We again obtain the bound
\begin{equation}\label{e:norm2}\norm{A_\ell A^*_m}\leq 
 {\frac{1}{(2\pi\hbar)^{r}}}
 \frac{ \sup_{\xi\in\widetilde\Omega_2}|(\det \nabla  p_n\circ\ldots\circ p_1(\xi))|}{
\inf_{\xi\in\widetilde\Omega_2}|(\det \nabla\tilde    p_n\circ\ldots\circ \tilde p_1(\xi))|} ,
\end{equation}
and $\norm{A_\ell A^*_m}=0$ if $\norm{\ell-m}>C$. Estimates \eqref{e:norm1} and \eqref{e:norm2}, combined with the Cotlar-Stein lemma, yield Theorem \ref{t:mainnorm2}. The last statement of the theorem comes from Remark \ref{r:unif}.

\section{Examples \label{s:exam}}
We now give examples of application of Theorems \ref{t:mainnorm} and \ref{t:mainnorm2}. 
These results are needed in \cite{AS} and \cite{AR10}.

Let $\bY$ be a $d$-dimensional $\cC^\infty$ manifold. The cotangent bundle $T^*\bY$ is endowed with its canonical symplectic form, denoted $\omega$. Let $H:T^* \bY\To\IR$ be a smooth function (hamiltonian), and let $\Phi_H^t:T^*\bY\To T^*\bY$ be the corresponding hamiltonian flow. We assume for simplicity that
$(\Phi_H^t)$ is complete. We fix a time step $\tau>0$, arbitrary. Before specifying the operators $\hat P_n$ to which we will apply the previous results, we have to make several assumptions concerning the underlying geometric situation.
 
We assume that we have a smooth foliation $\cF$ of $T^*\bY$ by lagrangian leaves (in the sequel we shall simply speak about a ``lagrangian foliation''), and will denote $\cF^{(n)}=\Phi^{n\tau}_H(\cF)$. Let $O\subset T^*\bY$ be an open, relatively compact subset of $T^* \bY$; we assume that we have a finite open covering of $O$, $O\subset O_1\cup O_2\cup\ldots\cup O_K$, with the following properties~:
\begin{itemize}
\item[]for any $n>0$ and any sequence $(\alpha_0, \alpha_1,\ldots, \alpha_{n-1})\in \{1,\ldots, K\}^n$
such that $O_{\alpha_0}\cap \Phi^{-\tau}_H(O_{\alpha_1})\cap\ldots\cap \Phi^{-(n-1)\tau}_H(O_{\alpha_{n-1}})\not=\emptyset$, we can find for all $k\leq n-1$ 
a smooth symplectic coordinate chart $\Psi_k:(O_{\alpha_k}, \omega)\To (\IR^{2d}, \omega_o)$ which maps $O_{\alpha_k}$ to a ball in $\IR^{2d}$, and the foliation $\cF^{(k)}\rceil_{O_{\alpha_k}}$ to the horizontal foliation of that ball; the collection of coordinate charts $(\Psi_k)$ may depend on the sequence $(\alpha_0, \alpha_1,\ldots, \alpha_{n-1})$, however, it can be chosen so that all the derivatives of $\Psi_k$ and $\Psi_k^{-1}$ are bounded, independently of $n$, $(\alpha_0, \alpha_1,\ldots, \alpha_{n-1})$, and $k$.

\end{itemize}

We can now give some details about the operators $\hat P_k$ to which we shall apply the main results. We fix a family $\hat\chi_1,\ldots, \hat \chi_K$ of $\hbar$-pseudodifferential operators (see the appendix), such that the full symbol of $\hat \chi_k$ is compactly supported inside $O_k$. We also assume that its principal symbol $\chi_k$ (which is a smooth function on $T^*\bY$) satisfies $\norm{\chi_k}_{\cC^0}\leq 1$.

Let $\hat H$ be a self-adjoint $\hbar$-pseudodifferential operator with principal symbol $H$. 

Fix, finally, a sequence $(\alpha_0, \alpha_2,\ldots, \alpha_{n-1})\in\{1,\ldots, K\}^n$.
We shall use Theorems \ref{t:mainnorm} and \ref{t:mainnorm2} to estimate the
norm of the product $\prod_{k=0}^{n-1}\hat \chi_{\alpha_{k+1}}e^{-\frac{i\tau \hbar \hat H}\hbar}\hat \chi_{\alpha_{k}}.$ The operator $\hat P_k$ will be $\hat \chi_{\alpha_{k+1}}e^{-\frac{i\tau \hbar \hat H}\hbar}\hat \chi_{\alpha_{k}}$, read in an adapted coordinate system~:

Once the sequence $(\alpha_0, \alpha_2,\ldots, \alpha_{n-1})$ is fixed, we consider the family of coordinates $\Psi_k$ described in our assumptions.
We fix a collection of Fourier integral operators $U_k~: L^2(\bY)\To L^2(\IR^d)$, associated with the canonical transformation $\Psi_k$ ($k=0,\ldots n-1$), and such that the pseudodifferential operator
$U^*_kU_k$ satisfies $U^*_kU_k\hat\chi_k=\hat\chi_k+\cO(\hbar^\infty)$ and $ \hat\chi_k=\hat\chi_k U^*_k U_k+\cO(\hbar^\infty)$(where the $\cO$ is to be understood in the $L^2(\bY)$-operator norm). Note that the operators $U_k$ depend on the sequence $(\alpha_0, \alpha_2,\ldots, \alpha_{n-1})$, but (in the geometric situation described above) we can assume that their symbols have derivatives of all orders bounded independently of $(\alpha_0, \alpha_2,\ldots, \alpha_{n-1})$ and of $k$. We take $\hat P_k=U_{k+1}\hat \chi_{\alpha_{k+1}}e^{-\frac{i\tau \hbar H}\hbar}\hat \chi_{\alpha_{k}}U_k^*$. It is a Fourier integral operator $L^2(\IR^d)\To L^2(\IR^d)$, associated with the canonical transformation $\kappa_k=\Psi_{k+1}\Phi_H^\tau\Psi_k^{-1}$, which by construction preserves the horizontal foliation. These operators also satisfy all assumptions (H), hence we can apply to them Theorems \ref{t:mainnorm} and \ref{t:mainnorm2}.

We give a concrete example of application, used in \cite{AS} and \cite{AR10}. 
Let $G$ denote a non-compact connected simple Lie
group with finite center.
We choose a Cartan involution $\Theta$ for $G$, and let $K < G$ be the
$\Theta$-fixed maximal compact subgroup. Let $\fG = Lie(G)$, and let $\theta$
denote the differential of $\Theta$, giving the Cartan
decomposition $\fG= \fK\oplus \fP$ with $\fK = Lie(K)$. 
Let $\bS = G/K$ be the associated symmetric space.
For a lattice $\Gamma < G$ we write $\bX =\Gamma\backslash G$ and
$\bY = \Gamma\backslash G/K$, the latter being a
locally symmetric space of non-positive curvature.

Fix now a maximal abelian subalgebra $\fA\subset\fP$. The dimension of $\fA$ is called the real rank of $G$, and will be denoted by $r$ in the sequel. We denote by $\fA^*$ the real dual of
$\fA$.
Let $\fG_\alpha  = \{X\in\fG, \forall H \in \fA : ad(H)X =  \alpha(H)X \}$,
$\Delta=t\Delta(\fA: \fG)=  \{  \alpha\in\fA^*  \setminus \{0\},  \fG_\alpha\not = \{0\}\}$ and call the latter the (restricted)
roots of $\fG$ with respect to $\fA$. For $\alpha\in\Delta$, we denote by $m_\alpha$ the dimension of  $\fG_\alpha$. The subalgebra $\fG_0$ is  $\theta$-invariant, and hence $\fG_0=
(\fG_0\cap\fP)\oplus (\fG_0\cap\fK)$. By the maximality of $\fA$ in $\fP$, we must then have $\fG_0 = \fA\oplus \fM$
where $\fM = Z_{\fK}(\fA)$, the centralizer of $\fA$ in $\fK$.

A subset $\Pi\subset\Delta(\fA: \fG)$ will be called a system of simple roots
if every root can be uniquely expressed as an integral
combination of elements of $\Pi$ with either all coefficients non-negative or
all coefficients non-positive.
 Fixing a simple system $\Pi$ we get a notion of positivity.
We will denote by $\Delta^+$ the set of positive roots,
by $\Delta^-=-\Delta^+$ the set of negative roots.
For $\fN = \oplus_{\alpha >0}\fG_\alpha$ and
$\lienb = \Theta\fN = \oplus_{\alpha <0}\fG_\alpha$ we have
$\fG = \lien \oplus \fA \oplus \fM \oplus \lienb$.
 
Let $N, A < G$ be the connected subgroups corresponding to the subalgebras
$\fN, \fA \subset\fG$ respectively, and let $M = Z_K(\fA)$.
Then $\fM = Lie(M)$, though $M$ is not
necessarily connected.

On $T^*\bS$, consider the algebra $\cH$ of smooth $G$-invariant
hamiltonians, that are polynomial in the fibers of the projection
$T^*\bS\To \bS$. The structure theory of semisimple Lie algebras shows that $\cH$ is
isomorphic to a polynomial ring in $r$ generators.  Moreover, the elements
of $\cH$ commute under the Poisson bracket. Thus, we have on $T^*\bS$ a
family of $r$ independent commuting Hamiltonian flows $H_1,..., H_r$.
Since all these flows are $G$-equivariant, they descend
to the quotient $T^*\bY$.

We apply the discussion above to $\bY$ and $H\in\cH$. We assume that $O$ is such that the differentials $(dH_1, \ldots, dH_r)$ are everywhere independent on $O$.  It is known that any given regular
common energy layer $\{H_1=E_1, \ldots, H_r=E_r\}\subset T^*\bY$ may naturally be identified (in a $G$-equivariant way) with $G/M$.
We thus have an equivariant map $O\To \IR^r\times G/M$ which is a diffeomorphism onto its image. 
In all that follows, we identify $O$ with an open subset of $\IR^r\times G/M$.
Under this identification, the action of $\Phi_H^t$ is transported to
$$(E_1, \ldots,E_r, \rho M)\mapsto (E_1, \ldots, E_r, \rho e^{t a_{E_1,\cdots, E_r}}M),$$
where $a_{E_1, \ldots, E_r}\in \fA$ depends smoothly on  $E_1, \ldots, E_r$, and linearly on $H$ -- see \cite{Hil, AS} for more explanations.
The foliation $\cF$ is invariant under $\Phi_H^t$, and can be described as follows~:
the leaf of $(E_1, \ldots, E_r, \rho M)\in \IR^r\times G/M$ is  $(E_1, \ldots, E_r)\times\{\rho a\bar n M, a\in A, \bar n\in\bar N\}.$

We assume that each $O_k$ is small enough so that, for any given $(E_1, \ldots, E_r, \rho M)\in O_k$, the map
\begin{eqnarray*}\IR^r \times \fN\times \fA\times \bar\fN&\To&\IR^r \times G/M\\
(\varepsilon_1, \ldots, \varepsilon_r, X, Y, Z)&\mapsto&( \varepsilon_1, \ldots, \varepsilon_r,\rho e^X e^Y e^Z M)
\end{eqnarray*}
is a local diffeomorphism from a neighbourhood of $(E_1, \ldots, E_r, 0,0,0)$ onto $O_k$. In such coordinates, the leaves of the foliation $\cF$ are then given by the equations $(\varepsilon_1, \ldots, \varepsilon_r, X)=cst$.

Let $(\alpha_0,\ldots, \alpha_{n-1})$ be such that $O_{\alpha_0}\cap \Phi^{-\tau}_H(O_{\alpha_1})\cap\ldots\cap \Phi^{-(n-1)\tau}_H(O_{\alpha_{n-1}})\not=\emptyset$.
We can then take $\rho=\rho_k$ and $(E_1, \ldots, E_r)$ such that $(E_1, \ldots, E_r, \rho_k M)\in O_{\alpha_k}$ and $(E_1, \ldots, E_r, \rho_{k+1} M)=\Phi_H^\tau(E_1, \ldots, E_r, \rho_k M).$
As explained in the previous paragraph, fixing $\rho_k$ allows to identify $O_{\alpha_k}$ with a subset of $\IR^r \times \fN\times \fA\times \bar\fN$; and we denote $(\varepsilon_1, \ldots, \varepsilon_r, X, Y, Z)$ these coordinates.

Denote by $d$ the dimension of $ \bS$ (note that $d=r+\dim \fN=r+\sum_{\alpha\in\Delta^+}m_\alpha$). By the Darboux-Lie theorem (see \cite{Vais}), we can find some coordinate system $\Psi_k=(x^k_1,\ldots, x^k_d, \xi^k_1,\ldots, \xi^k_d):O_{\alpha_k}\To \IR^d$ mapping $\omega$ to $\omega_o$, and such that $ (\xi^k_1,\ldots, \xi^k_d)=(\varepsilon_1,\cdots, \varepsilon_r, X)$. The canonical transformations $\kappa_k=\Psi_{k+1}\Phi_H^\tau\Psi_k^{-1}$
preserve the horizontal foliation of $\IR^d$, hence they are of the form $\kappa_k~: (x, \xi)\mapsto (x', \xi'=p_k(\xi)).$ It turns out, in this particular case, that the maps $p_k$ are all the same, and of the particular form
$$p_k(\varepsilon_1,\cdots, \varepsilon_r, X)=(\varepsilon_1, \ldots, \varepsilon_r, Ad(e^{\tau a_{\varepsilon_1,\cdots, \varepsilon_r}}).X),$$
where  $a_{\varepsilon_1, \ldots, \varepsilon_r}\in \fA$ depends smoothly on  $\varepsilon_1, \ldots, \varepsilon_r$.
The linear maps $Ad(e^{\tau a_{\varepsilon_1, \ldots, \varepsilon_r}})$ acting on $\fN$ are all simultaneously diagonalizable, the eigenspaces being the root spaces $\fG_\alpha$, with eigenvalue $e^{\tau\alpha(a_{\varepsilon_1, \ldots, \varepsilon_r})}$. We are thus in a case of application of Theorem \ref{t:mainnorm2} provided that $\alpha(a_{\varepsilon_1, \ldots, \varepsilon_r})\leq 0$ for all $\alpha\in \Delta^+$. If we fix $J\subset\Delta^+$ arbitrarily,  the map $\kappa_k$ preserves the coisotropic foliation by the leaves $(\varepsilon_1, \ldots, \varepsilon_r, X_J)=cst$ (where any $X\in\fN$ is decomposed into $X=\sum_{\alpha\in \Delta^+}X_\alpha$, $X_\alpha\in\fG_\alpha$, and $X_J$ is defined by $X_J=\sum_{\alpha\in J}X_\alpha$).

\begin{cor} Assume that $H$ and $O$ are such that $\alpha(a_{\varepsilon_1, \ldots, \varepsilon_r})\leq \delta\leq 0$ for all $(\varepsilon_1, \ldots, \varepsilon_r, \rho M)\in O$ and all $\alpha\in\Delta^+$.
Fix a subset $J\subset\Delta^+$ of the sets of roots.

Fix $\cK>0$ arbitrary. Then, for $n=\cK |\log \hbar|$, and for every sequence $(\alpha_0, \ldots, \alpha_{n-1})$,
$$\norm{\prod_{k=0}^{n-1}\hat \chi_{\alpha_{k+1}}e^{-\frac{i\tau \hbar H}\hbar}\hat \chi_{\alpha_{k}}}\leq
\frac{1}{(2\pi\hbar)^{(d-r-\sum_{\alpha \in J}m_\alpha)/2}}\prod_{\alpha\in\Delta^+\setminus J}e^{\frac{n\delta}2}
$$
for $\hbar>0$ small enough.
\end{cor}
 \begin{rem} In the situation of \cite{AS}, we actually do not have $\alpha(a_{\varepsilon_1, \ldots, \varepsilon_r})\leq \delta\leq 0$, but $\alpha(a_{\varepsilon_1, \ldots, \varepsilon_r})\leq \delta$, where $\delta>0$ can be made arbitrarily small by conveniently choosing the set $O$. Once $\cK$ is given, we can choose $\delta$ small enough (and $O$) so that the proof of Section \ref{s:proof} still works for $n=\cK |\log \hbar|$.
 \end{rem}
 In the special case $G=SO_o(d, 1)$, $\bY$ is a hyperbolic manifold of dimension $d$. We have $r=1$, and $\cH$ is generated by the laplacian $\lap$. Taking $H=-\frac{\lap}2$, $J=\emptyset$, and keeping the same notations as before, we obtain in this case
 $$\norm{\prod_{k=0}^{n-1}\hat \chi_{\alpha_{k+1}}e^{-\frac{i\tau \hbar H}\hbar}\hat \chi_{\alpha_{k}}}\leq
\frac{1}{(2\pi\hbar)^{(d-1)/2}} e^{-n\tau\frac{(d-1)}2 (1-\eta)}
$$
if we assume that the symbols of the pseudodifferential operators $\hat\chi_k$ are all supported in
$\{\norm{\xi}\in [1-\eta, 1+\eta]\}$ for some small $\eta>0$.
The result proved in \cite{An, AN07-1}, which was only based on the idea of Theorem \ref{t:mainnorm}, was
$$\norm{\prod_{k=0}^{n-1}\hat \chi_{\alpha_{k+1}}e^{-\frac{i\tau \hbar H}\hbar}\hat \chi_{\alpha_{k}}}\leq
\frac{1}{(2\pi\hbar)^{d/2}} e^{-n\tau\frac{(d-1)}2(1-\eta)}.
$$
We see that Theorem \ref{t:mainnorm2} allows to improve the prefactor $\frac{1}{(2\pi\hbar)^{d/2}}$ to $\frac{1}{(2\pi\hbar)^{(d-1)/2}} $, as needed in \cite{AR10}.

\begin{rem} Versions of the hyperbolic dispersion estimate have also been proved for more
general uniformly hyperbolic dynamical systems \cite{An, AN07-1, NZ09, GR09}, and even
for certain non-uniformly hyperbolic systems \cite{GR10}. We refer the reader to \cite{N10} for an expository paper. It is not clear to me whether the new presentation (and improvement) introduced here can be used for those systems. Indeed,
there is in general no {\em smooth} lagrangian foliations preserved by the hamiltonian flow, and so one cannot hope that  the symplectic changes of coordinates $\Psi_k$ used above will have uniformly bounded derivatives. On the other hand, control of high order derivatives is crucial when one applies the techniques of semiclassical
analysis (method of stationary phase, integration by parts,...) It is a drawback of semiclassical analysis that it cannot deal with symplectic transformations of low regularity~: I don't know if this obstacle can be overcome.
\end{rem}

 \section{A few definitions\label{s:def}} 
 For the reader's convenience, we clarify the terminology used in this paper. This is of course
 not an exhaustive tutorial on semiclassical analysis~: for this we refer the reader to \cite{EZ}.
 
 Let $M$ and $N$ be two smooth manifolds of the same dimension $d$. Their cotangent bundles $T^*M$ and $T^*N$ are respectively equipped with  the canonical symplectic forms $\omega_M$ and $\omega_N$. Let $\kappa$ be a diffeomorphism from an open relatively compact subset $O_M\subset T^*M$ onto an open subset $O_N\subset T^*N$, sending $\omega_M$ to $\omega_N$ (such a $\kappa$ is called a symplectic diffeomorphism, of a canonical transformation).
 
 In this paper, we say that an operator\footnote{More properly, a family of operators depending on $\hbar>0$} $\hat P=\hat P_\hbar : L^2(M)\To L^2(N)$ is a (semiclassical) Fourier integral operator associated with $\kappa$, if it is a finite sum of operators of the form
 $$\hat Q f(x')=\frac1{(2\pi\hbar)^{\frac{d+m}2}}\int_{x\in\IR^d, \theta\in\IR^m} e^{\frac{i S(x, x', \theta)}\hbar}a(x, x', \theta, \hbar)f(x)dx d\theta, $$ where
 \begin{itemize}
\item $m\geq 0$ is an integer and $S$ is a smooth function on $M\times N\times \IR^m$.  
\item For a given $\hbar>0$, the function $(x, x', \theta)\mapsto a(x, x', \theta, \hbar)$ is of class $\cC^\infty$ and has compact support, independent of $\hbar$;
\item When $\hbar \To 0$, $a(x, x', \theta, \hbar)$ has an asymptotic expansion
$$a(x, x', \theta, \hbar)\sim   \sum_{k=0}^{\infty}\hbar^k a_k(x, x', \theta),$$
valid up to any order and in all the $\cC^\ell$ norms.
\item for any $(x, x', \theta)$ in the support of $a$, and for any $(x, \xi)\in T^*_xM$, $(x', \xi')\in T^*_{x'}N$, we have $$\big[\xi=-\partial_x S_n(x, x', \theta), \;\xi'=\partial_{x'} S_n(x, x', \theta),\;
\partial_\theta S_n(x, x', \theta)=0\big]\Longrightarrow\big[(x, \xi)\in O_M, ( x', \xi')=\kappa(x, \xi)\big].$$
\end{itemize}
If $\hat P$ is a Fourier integral operator associated with $\kappa$ and $\hat P'$ is a Fourier integral operator associated with $\kappa'$, then $\hat P'\circ\hat P$ is a Fourier integral operator associated with $\kappa'\circ\kappa.$ Also, the adjoint $\hat P^*$ is a Fourier integral operator associated with $\kappa^{-1}.$

In the case $M=N$, a pseudodifferential operator is a Fourier integral operator associated with $\kappa=$ the identity.
 
 If $\hat P$ is a pseudodifferential operator on $M=\IR^d$, we define its full symbol as the function $a_\hbar(x, \xi)$, defined on $\IR^d\times\IR^d$, by
 $$a_\hbar(x_0, \xi_0)=e^{\frac{-i\la\xi_0, x_0\ra}\hbar}\hat P(e^{\frac{i\la\xi_0, \bullet\ra}\hbar})\rceil_{x_0}.$$ By the stationary phase method, the function $a_\hbar$ can be shown to have an asymptotic expansion, valid in all $\cC^k$-norms on $\IR^d$,
 $$a_\hbar(x, \xi)\sim\sum_{j\in\IN} \hbar^j a_j(x, \xi).$$
 The term $a_0$ is called the principal symbol.
 We say that the full symbol of $\hat P$ vanishes on $\Omega\subset \IR^d\times\IR^d$ if all the $a_j$ vanish on $\Omega$.


\begin{thebibliography}{a}

 
\bibitem{An} N.~Anantharaman, {\em Entropy and the localization of eigenfunctions},  Ann. of Math. (2)  168  (2008),  no. 2, 435--475.

\bibitem{ICM} N.~Anantharaman, {\em A hyperbolic dispersion estimate, with applications to the linear Schr\"odinger equation}, to appear in Proceedings of the ICM 2010.


\bibitem{AN07-1}N.~Anantharaman, S.~Nonnenmacher, \emph{Half-delocalization of eigenfunctions for the laplacian on an Anosov manifold}, Festival Yves Colin de Verdi\`ere.  Ann. Inst. Fourier (Grenoble)  57  (2007),  no. 7, 2465--2523. 

\bibitem{AKN07}
N.~Anantharaman, H.~Koch and S.~Nonnenmacher, \emph{Entropy of eigenfunctions}, New Trends of Mathematical Physics, selected contributions of the 15th International Congress on Mathematical Physics, Springer (2009), 1--22.

\bibitem{AR10} N.~Anantharaman, G.~Rivi\`ere,  \emph{Dispersion and controllability for the Schr\"odinger equation on negatively curved manifolds}, preprint 2010

\bibitem{AS} N.~Anantharaman, L.~Silberman, \emph{Asymptotic distribution of eigenfunctions
on locally symmetric spaces}, work in progress.

\bibitem{BGH} N.~Burq, C.~Guillarmou, A.~Hassell, \emph{Strichartz estimates without loss on manifolds with hyperbolic trapped geodesics}, preprint, arXiv:0907.3545.

\bibitem{Ch1} H.~Christianson, \emph{Semiclassical non-concentration near hyperbolic orbits}, J. Funct. Anal. 262 (2007),
145--195; ibid, \emph{Dispersive estimates for manifolds with one trapped orbit,} Comm. PDE 33 (2008), 1147--1174.


\bibitem{Ch2} H.~Christianson, \emph{Cutoff resolvent estimates and the semilinear Schr\"odinger equation}, Proc. AMS 136 (2008),
3513--3520.

\bibitem{Ch3} H.~Christianson, \emph{Applications to cut-off resolvent estimates to the wave equation}, Math. Res. Lett. Vol. 16 (2009), no. 4, 577--590.

\bibitem{Dat} K.~Datchev, \emph{Local smoothing for scattering manifolds with hyperbolic trapped sets}, Comm. Math. Phys. 286, no. 3, 837--850. 

\bibitem{Doi} S.--I. Doi, \emph{Smoothing effects of Schr\"odinger evolution groups on Riemannian manifolds}, Duke Math. J.
82( 1996), 679--706.


\bibitem{EZ} L.C.~Evans, M.~Zworski \emph{Lectures on semiclassical analysis} (version 0.3) avalaible at http://math.berkeley.edu/~zworski/semiclassical.pdf (2003)


 \bibitem{Hil} J.~Hilgert 
{\it An ergodic Arnold-Liouville theorem for locally symmetric spaces.} Twenty years of Bialowieza: a mathematical anthology, 163--184,
World Sci. Monogr. Ser. Math., 8, World Sci. Publ., Hackensack, NJ, 2005.

\bibitem{N10} S.~Nonnenmacher, \emph{Entropy of chaotic eigenstates}, Notes of the minicourse given at the workshop ``Spectrum and dynamics'', Centre de Recherches Math\'ematiques, Montr\'eal, April 2008, arXiv:1004.4964.
  
\bibitem{NZ09} S.~Nonnenmacher, M.~Zworski, 
\emph{Quantum decay rates in chaotic scattering}, Acta Mathematica, Volume 203, Number 2, December 2009, 149--233.

\bibitem{NZ09-1} S.~Nonnenmacher, M.~Zworski, 
\emph{Semiclassical resolvent estimates in chaotic scattering}, Appl. Math. Res. Express (2009) 74--86.

 \bibitem{GR09} G.~Rivi\`ere, \emph{Entropy of semiclassical measures in dimension 2}, to appear in Duke Math. J.

\bibitem{GR10} G.~Rivi\`ere, \emph{Entropy of semiclassical measures for nonpositively curved surfaces}, preprint.

\bibitem{Vais} I.~Vaisman, \emph{Basics of Lagrangian foliations},  Publ. Mat.  33  (1989),  no. 3, 559--575.

 



 \end{thebibliography}
\end{document}